\documentclass{article}[12pt]
\usepackage[cp1251]{inputenc}
\usepackage{amsmath,amsfonts,amssymb,amsthm,mathrsfs,bm,cleveref,color}

\usepackage{graphicx}

\usepackage{caption}
\usepackage{csquotes}
\usepackage{cite}

\renewcommand{\Im}{\mathop{\mathrm{Im}}\nolimits}

\begin{document}

\noindent
\textbf{Physical optics}
\vskip .2cm

\noindent
{\large\textbf{UNIDIRECTIONAL PULSES:\\ RELATIVELY UNDISTORTED QUASI-SPHERICAL WAVES, FOURIER-BESSEL INTEGRALS,\\
 AND  PLANE-WAVES   DECOMPOSITIONS}
}

\vskip .3cm
\noindent \textbf{Alexandr B. Plachenov}$^1$,
 \textbf{Aleksei P. Kiselev}$^{2,3*}$
\vskip .5cm
{\scriptsize
\noindent$^1$ \emph{MIREA -- Russian Technological University, Moscow, 119454 Russia}

\noindent$^{2}$ \emph{St. Petersburg Department of V.A. Steklov Mathematical Institute of the Russian Academy of Sciences, St. Petersburg, 191023 Russia;
\\
$^{3}$Institute for Problems in Mechanical Engineering of the Russian Academy of Sciences,
St. Petersburg, 199178 Russia
}}

\noindent \emph{*e-mail: aleksei.kiselev@gmail.com}
\vskip .3cm
\noindent {\small

\vskip .3cm
\noindent \textbf{Abstract}:
A theoretical description of a class of unidirectional axisymmetric localized pulses, is given. The equivalence of their representations in the form of relatively undistorted quasi-spherical waves, in the form of Fourier-Bessel integrals and in the form of a superposition of plane waves with wave vectors having positive projections on a given direction is established.

\vskip .3cm
}

\section{Introduction}In recent years, there has been an increased interest in localized solutions of the wave equation
\begin{equation}\label{we}
\frac{\partial^2u}{\partial x^2}+\frac{\partial^2u}{\partial y^2}+\frac{\partial^2u}{\partial z^2}-\frac{1}{c^2}\frac{\partial^2u}{\partial t^2}=0,
\end{equation}(with $x$, $y$, and $z$ Cartesian coordinates,  $t$  time, and $c > 0$ wave speed assumed constant),
having the property of unidirectionality
 \cite{Zamboni,PaAl,Lekner,Lekner_book,SPK1,SPK2,Birula,Pla20,BeSa,Lekner23}.
Real and imaginary parts of such solutions can be used as components of the Hertz's vector in the construction of unidirectional electromagnetic pulses.

One of the formulations of unidirectionality (see, for example, \cite{Pla20}) consists in the requirement that only homogeneous plane waves traveling in directions forming an angle with a certain chosen direction not exceeding $\frac{\pi}{2}$, are present in the decomposition of the solution in plane waves. This property expresses the requirement, natural from a physical point of view, that the mathematical model of the pulse describe its propagation strictly from the source. Unidirectional pulses are sometimes called causal \cite{Lekner,Lekner_book}.
It is noteworthy that the unidirectionality understood in the above sense does not exclude the possibility
that in some spatiotemporal regions the projection of the energy flow vector
on the chosen direction may turn out to be negative, see \cite{Birula,Berry,BeSa}.

In what follows, the chosen direction of propagation will be the direction of the $z$ axis.
Accordingly, we will call the solution unidirectional if the $z$ projections of the wave vectors of its plane-wave constituents are non-negative. Trivial examples of unidirectional solutions are a plane wave and a finite combination of plane waves that are not spatially localized. In this paper, we address solely unidirectional pulses,
localized with respect to all spacial coordinates at any fixed instant in time.

The first results on the construction of unidirectional pulses were based on the consideration of axisymmetric solutions of the equation \eqref{we} in the form of Fourier-Bessel integrals
\begin{equation}\label{FouBess}
 u= u(\rho,z,t)=\int_0^\infty d\omega e^{i\omega t} \int_0^{{\omega}/{c}}dk_z A(k_z,\omega)e^{-ik_zz}
  J_0(\rho\sqrt{{\omega^2}/{c^2}-k_z^2}),
\end{equation}
where $\rho=\sqrt{x^2+y^2}$,
     with fairly arbitrary weight functions $A$ \cite{Zamboni,PaAl,Lekner}. A proper choice of such a weight allowed  to find several
solutions expressed in terms of elementary functions. The simplest localized unidirectional solution, however, was found differently and was based on a lucky trick \cite{SPK1,SPK2}, which used a partial fraction decomposition of the well-known splash pulse (see, for example, \cite{Ziolk,Feng,2D}). This solution is
\begin{equation}\label{IS}
      u= \frac{1}{S(S-z_*)}\,,
    \end{equation}
where
\begin{equation}\label{S}
S=S(t,\textbf{R})=\sqrt{c^2 t_*^2-\rho^2},\end{equation}
with $\textbf{R}$ denoting the position vector of the observation point, and
\begin{equation}\label{**}
z_*=z+i\zeta, \,\,t_*=t+i\tau\,.
\end{equation}
Here, $\zeta$ and $\tau>0$ are free parameters which are assumed real. The square root branch in \eqref{S} is chosen so that $\left. S\right|_{x=y=0}=ct_*$,
in which case
for arbitrary values of $t,\textbf{R}$ the inequality $\Im S\ge c\tau$  \cite{SPK2} holds.
  The solution \eqref{IS} is non-singular under the condition  $\zeta<c\tau$, and then its energy is finite \cite{SPK2,Pla20}. We also note that
with proper choice of free parameters $\zeta$ and $\tau$, this solution can model pancake-shaped, ball-shaped and needle-shaped focused pulses.
   The proof of its unidirectionality was presented in \cite{Pla20,BeSa}.

On the other hand, Besieris and Saari \cite{BeSa} (see also \cite{PaAl,Lekner23} noted that a special class of \emph{relatively undistorted waves} is important in the description of unidirectional wave propagation. This is the name of the solutions of the equation \eqref{we} of the form
\begin{equation}\label{}
u=gf(\theta)\,,
\end{equation}
where the \emph{amplitude} $g=g(x,y,z,t)$ and the \emph{phase} $\theta=\theta(x,y,z,t)$ functions are fixed, and \emph{waveform} $f$ is an \emph{arbitrary} function \cite{CH,K07}.
If the phase function $\theta$ is complex-valued, then the waveform $f$ must be analytical in the range of values $\theta$ \cite{KPD}.

The class of axisymmetric solutions we are interested in has
the form
    \begin{equation}\label{RNDW}
    u= \frac{f(\theta)}{S}\,,
    \end{equation}
    where
    \begin{equation}\label{thetaqsph}
    \theta=S-z-ib,
    \end{equation}
    $b=c\tau>0$, with an arbitrary waveform $f$ analytical in the upper half-plane $\mathbb{C}^+$, see \cite{SPK2}.
   The imaginary shift $ib$ in \eqref{thetaqsph} is introduced in order to have the upper half-plane as the range of values of the phase $\theta$.

In order for solutions of the form \eqref{RNDW} to describe a localized wave,
we will require a fairly fast decrease (no slower than $1/|\theta|$) of the function $|f(\theta)|$ when $|\theta|\to \infty$.

Let us explain how one can easily come to such solutions. Let us consider a relatively undistorted wave corresponding to spherical waves
   \cite{CH},
 \begin{equation}\label{sph}
      u= \frac{f(R-ct)}{R},
    \end{equation}
$R=\sqrt{x^2+y^2+z^2}$. The equation \eqref{we} is invariant under replacement
    \begin{equation}\label{ctzb}
    z \mapsto i(ct+ib)\,, \qquad ct \mapsto i(z+ib)\,.
    \end{equation}
where $b$ is a free real constant. To match with \eqref{IS}-\eqref{**}, we take $b=c\tau>0$. Under the transformation \eqref{ctzb},  the phase function $R-ct$ in \eqref{sph}   goes into the expression \eqref{thetaqsph}, which
we will call it \emph{quasi-spherical phase}. The imaginary part of \eqref{thetaqsph} is non-negative for all values of
     spatial and temporal variables. Redefining the waveform using the rule $f(i\theta)/i \mapsto f(\theta)$,
     we obtain a class of relatively undistorted waves of the form \eqref{RNDW},
     which we will call \emph{quasi-spherical}.

Note that the expression \eqref{IS} is a special case of \eqref{RNDW} when
$f(\theta)=-\frac{1}{\theta+i(b-\zeta)}$.
   A number of more complex, although expressed through elementary functions, solutions from this class were found in \cite{Zamboni,Lekner,Lekner23}.
   A solution close to that found in \cite{SPK1,SPK2}, but less general -- with
  $\zeta=0$ -- is presented in \cite{Birula}. The article \cite{BeSa} provides an overview of such solutions.

 An example of a unidirectional solution for a time-harmonic regime is given in \cite{KP19}.
  This solution has an asymptotic behavior corresponding to a Gaussian beam with arbitrary astigmatism.

In this note, we establish a relation between solutions described by  \eqref{RNDW} and \eqref{FouBess}. Further,
we represent these solutions in the form of a superposition of plane waves. Our approach is based on a technique for studying localized solutions that goes back to  Blagoveshchenskii \cite{Blagoveshchenskii} and Moses--Prosser \cite{MP}. It rests upon  formulas expressing the solution through its asymptotic behavior in the far zone at large time. This technique  turns  out to be convenient, in particular, for  calculating such characteristics of a localized pulse as energy, momentum and orbital angular momentum \cite{Pla20,PlCPK23}.
\section{The Blagoveshchenskii --- Moses --- Prosser approach and the unidirectionality of quasi-spherical waves}

\subsection{Blagoveshchenskii --- Moses --- Prosser approach}
We denote the position vector by  ${\bf R}=x\textbf{\emph{i}}+y\textbf{\emph{j}}+z\textbf{\emph{k}}$, where
$\emph{\textbf{i}}$, $\emph{\textbf{j}}$ and $\emph{\textbf{k}}$ are unit vectors along the coordinate axes. Let   ${\bf n}=\frac{\bf R}{R}$, $|{\bf n}|=1$ be the related to ${\bf R}$  unit vector, and $R$ = $|{\bf R}|=\sqrt{x^2+y^2+z^2}$ be the distance to the Cartesian origin.

 Consider an arbitrary smooth localized solution of the wave equation \eqref{we}, assuming that $R$ and $ct$ \emph{grow consistently}, i.e., that their difference
\begin{equation}\label{s}
  s = R - ct
\end{equation}
remains constant. Evidently,
$\textbf{R} = (ct + s)\textbf{n}$.

In the works of Blagoveshchenskii \cite{Blagoveshchenskii} and Moses--Prosser \cite{MP}, it was found that for any solution of the wave equation decreasing rapidly enough at $R\to\infty$, for any fixed $s$ and any direction $\textbf{n}$ there
exits the limit
\begin{equation}\label{F(s,n)}
F(s,{\bf n})=\lim_{t\to\infty} \left[ct \, u\left(t, (ct+s){\bf n}\right) \right].
\end{equation}
The limit \eqref{F(s,n)} characterizes, at large values of time, the amplitude
of the pulse in the direction of $\textbf{n}$.
For a unidirectional (along the $z$ axis) wave packet, obviously
${F(s,\textbf{n})}\equiv 0$ for all $\textbf{n}$ whose projections on the $z$ axis are negative, i.e. $\textbf{n}\cdot\textbf{\emph{k}}<0$,
where ${\bf n}\cdot \textbf{\emph{k}}$ is the scalar product of the vectors ${\bf n}$ and $\textbf{\emph{k}}$.

We will characterize the direction ${\bf n}$ by the angles $\chi$ and $\varphi$ of the spherical coordinate system with the polar axis $z$:
$${\bf n}=
\sin\chi\cos\varphi\, \textbf{\emph{i}} + \sin\chi\sin\varphi\, \textbf{\emph{j}}+ \cos\chi\,\textbf{\emph{k}},$$
$0\le\varphi< 2\pi$, $0\le\chi\le\pi$. The unidirectionality condition takes the form
\begin{equation}\label{<0}
F(s,{\bf n})\equiv 0 \,, \, \, \, \frac{\pi}{2}<\chi\leq\pi.
\end{equation}

The nontrivial result of Blagoveshchenskii -- Moses -- Prosser (see \cite{Blagoveshchenskii,MP,K07}) is that the solution of $u$ at any point of $\textbf{R}$ at any moment $t$ is representable
through the limit \eqref{F(s,n)} as follows:
\begin{equation}\label{plane}
u(t,{\bf R})= \frac{1}{2\pi} \iint_{|{\bf N}|=1} F'({\bf N}\cdot {\bf R}- ct,{\bf N})\, d^2{\bf N}\,,
\end{equation}
where the notation is introduced
$$F'(s,{\bf N})=\frac{\partial F(s,{\bf N})}{\partial s}.$$
The integration is carried out over the unit sphere $|{\bf N}|=1$, and $ d^2{\bf N}$ denotes the element of its surface area. In spherical coordinates
$\textbf{N}=\sin{\cal{X}}\cos\phi\, \textbf{\emph{i}} + \sin{\cal{X}}\sin\phi\, \textbf{ \emph{j}}+ \cos{\cal{X}}\, \textbf{\emph{k}},$
and the area element of the sphere takes the form
 $d^2{\bf N}=\sin{{\cal{X}}}\,d{\cal{X}}\, d\phi$.

The formula \eqref{plane} represents the solution $u$ in the form of a superposition of nonstationary plane waves.

\subsection{Quasi-spherical wave \eqref{RNDW} at large time and large distance}

Let us find the limit \eqref{F(s,n)} for the solution
 \eqref{RNDW}.
 If $\cos\chi\neq0$, then at $t\to+\infty$
\begin{gather*}S{=
S(t, (ct+s){\bf n})}=\sqrt{(ct+ib)^2-(ct+s)^2\sin^2{\chi}} \\
\approx ct|\cos\chi|+ \frac{ib-s\sin^2\chi}{|\cos\chi|}
\approx ct|\cos \chi|\,,
\end{gather*}
so that
\begin{gather*}
\theta=S-(ct+s)\cos\chi-ib\\
\approx ct(|\cos\chi|-\cos\chi) +
\frac{-s(\sin^2\chi+\cos\chi|\cos\chi|)+ib(1-|\cos\chi|)}{|\cos\chi|}\,.
\end{gather*}

For directions ${\bf n}$ making an obtuse angle with the $z$ axis,
   $\cos\chi<0$,
that is, $\chi>\frac{\pi}{2}$, we have
$u\approx {f(2ct|\cos\chi|)}/{(ct|\cos \chi|)}\,$,
    and since $f(2ct|\cos\chi|) \to 0$, from \eqref{F(s,n)}
     it follows that $$F(s, {\bf n})=0.$$
For directions ${\bf n}$ making an acute angle with the   $z$ axis,
$\cos\chi>0$,
$$u(t, (ct+s){\bf n})\approx \frac{1}{ct\cos\chi}f\left(\frac{-s+ib(1-\cos\chi)}{\cos\chi}\right)\,,$$
and \eqref{RNDW} implies
$$F(s, {\bf n})=\frac{1}{\cos\chi}f\left(\frac{-s+ib(1-\cos\chi)}{\cos\chi}\right)\,.$$

Finally, for $\chi=\frac{\pi}{2}$ a similar calculation gives the value
\begin{equation*}
F(s,{\bf n})=\lim_{t\to\infty} \left[\sqrt{\frac{ct}{2(ib-s)}} f\left(\sqrt{2ct(ib-s)} \right) \right]
= \frac{1}{2(ib-s)}\lim_{\theta\to\infty} \theta f(\theta).
\end{equation*}
    This limit is finite if $|f(\theta)|$ decreases not slower than $|\theta|^{-1}$, as we have assumed.
     Since the circle $\chi=\frac{\pi}{2}$ does not contribute to the integral \eqref{plane}, the value of $F(s,{\bf n})|_{\chi=\frac{\pi }{2}}$ can be replaced by zero and the result  written in the form
\begin{equation}\label{FSo}
F(s,{\bf n})=\frac{H(\cos\chi)}{\cos\chi}
f\left(\frac{-s+ib(1-\cos\chi)}{\cos\chi}\right)\,,
\end{equation}
with $H$ the Heaviside step function
\begin{equation}
	H(p) = \begin{cases}
		1, \, p > 0,\\
		0,\, p < 0.
	\end{cases}
\end{equation}

Thus, since for a quasi-spherical wave in the formula \eqref{plane} the integration occurs over the forward hemisphere ${\bf N}\cdot\emph{\textbf{k}}>0$,
the unidirectionality of the pulse \eqref{RNDW} is established.

\subsection{On angular divergence of quasi-spherical waves }
It should be noted that quasi-spherical solutions can describe pulses having not only a significant (as in the examples discussed in \cite{SPK1,SPK2,Birula}), but also a small angular divergence.
Strong angular localization requires a rapid decrease in the modulus of the function $f(\theta)$ with the growth of $\Im \theta$. This property is possessed, for example, by the waveform introduced by Lekner
    \cite{Lekner23}, having the form $$f(\theta)= \exp(iK\theta))/(\theta+ib),$$
    ($K$ is a real constant), for which the angular localization in the angle $\chi$ has
    Gaussian character.

    A number of examples of a waveform that provides Gaussian localization not only in angles, but also in a longitudinal variable can be found
     in the work of Kiselev and Perel \cite{PK2} (see also \cite{K07,PK1}), devoted to wave packets of a different nature.

\section{Integral representations of the quasi-spherical wave}

\subsection{Representation by superposition of non-stationary plane waves}
Differentiating the function \eqref{FSo} with respect to the first argument, we obtain
\begin{equation}\label{'}F'(s, {\bf n})
=-\frac{H(\cos\chi)}{\cos^2\chi}f'\left(\frac{{-s+ib(1-\cos\chi)}}{\cos\chi}\right)
\,,\,\chi\neq\frac{\pi}{2}.
\end{equation}
 Substituting \eqref{'} into \eqref{plane} and replacing $\textbf{n}$ with $\textbf{N}$ gives
\begin{gather} \nonumber%
u
=-\frac{1}{2\pi} \iint_{\Sigma_+} \frac{d^2{\bf N}}{\cos^2{\cal{X}}} f'\left( \frac{ct-{\bf N}\cdot {\bf R} +ib(1-\cos{\cal{X}})}{\cos{\cal{X}}}
\right)\\
=-\frac{1}{2\pi} \int_{-\pi}^{\pi}\,d\phi
\int_{0}^{\frac{\pi}{2}} \frac{\sin{\cal{X}}\,d{\cal{X}}}{\cos^2{\cal{X}}}
    f'\left( \frac{(ct+ib)-(z+ib)\cos{\cal{X}}-(x\cos\phi+y\sin\phi )\sin{\cal{X}}}{\cos{\cal{X}}}\right)
\,, \label{uF1}
\end{gather}
where $\Sigma_+$ denotes the forward semi-sphere $\{|{\bf N}|=1, \,{\cal{X}}<\frac{\pi}{2}\}$.

\subsection{Representation by a superposition of monochromatic plane waves }

Let us represent the waveform as follow
$$f(\theta)=\int_{0}^{\infty} \hat{f}(\kappa) \exp(i\kappa \theta) \,d\kappa\,,$$
this representation holds for fairly rapidly decreasing functions  $f$.
Integration is carried out along the positive semi-axis in view of the analyticity of $f$ in the upper half-plane. Then
$$f'(\theta)=i \int_{0}^{\infty} \hat{f}(\kappa) \exp(i\kappa \theta) \kappa \,d\kappa\,.$$
Therefore,
\begin{multline}\label{uk}
u
=-\frac{i}{2\pi} \int_{-\pi}^{\pi}\,d\phi
\int_{0}^{\frac{\pi}{2}} \frac{\sin{\cal{X}}\,d{\cal{X}}}{\cos^2{\cal{X}}} \\
\int_{0}^{\infty} \hat{f}(\kappa) \exp\left(i\kappa \frac{(ct+ib)-(z+ib)\cos{\cal{X}}-(x\cos\phi+y\sin\phi )\sin{\cal{X}}}{\cos{\cal{X}}}\right) \kappa \,d\kappa \\
=-\frac{i}{2\pi} \int_{-\pi}^{\pi}\,d\phi
\int_{0}^{\frac{\pi}{2}} {\sin{\cal{X}}\,d{\cal{X}}} \\
\int_{0}^{\infty} \hat{f}(k \cos{\cal{X}}) \exp\left[ik \left({(ct+ib)-(z+ib)\cos{\cal{X}}-(x\cos\phi+y\sin\phi )\sin{\cal{X}}}\right)\right] k \,dk
\,,
\end{multline}
   where the replacement  $\kappa=k \cos{\cal{X}}$ was made.
The right-hand side of \eqref{uk} can be understood as a volume integral  presented in spherical coordinates
 $(k,{\cal{X}},\phi)$, with  $k^2\,dk\,\sin{\cal{X}}\, d{\cal{X}}\, d\phi$ being the volume element and the integrand being
$$\hat{f}(k \cos{\cal{X}}) \frac{e^{ik[(ct+ib)-(z+ib)\cos{\cal{X}}-(x\cos\phi+y\sin\phi)\sin{\cal{X}}]}}{k}.$$
 The area of integration is the half-space  $0\le{\cal{X}}< \frac{\pi}{2}$.
Passing in \eqref{uk} to Cartesian coordinates   $$k_z=k\cos{\cal{X}}\,,\,\,
k_x=k\sin{\cal{X}}\cos\phi\,,\,\, k_y=k\sin{\cal{X}}\sin\phi\,,$$  we find
\begin{equation}\label{Fourier}
u
=-\frac{i}{2\pi} \int_{0}^{\infty} \hat{f}(k_z) dk_z \int_{-\infty}^{\infty} dk_x \int_{-\infty}^{\infty} dk_y
\frac{e^{i[k(ct+ib)-k_z(z+ib)-k_xx-k_yy]}}k\,,
\end{equation}
with $k=\sqrt{k_x^2+k_y^2+k_z^2}$.

Thus, we presented  the quai-spherical wave
 \eqref{RNDW} in the form of expansion in monochromatic plane waves.

   Continuing the integrand with zero to negative values of $k_z$, we obtain
\begin{equation}\label{Fourier3D}
u
=-\frac{i}{2\pi} \iiint_{\mathbb{R}^3}H(k_z) \hat{f}(k_z)
\frac{e^{i[k(ct+ib)-k_z(z+i\zeta)-k_xx-k_yy]}}k\,d^3{\bf k},
\end{equation}
with $d^3{\bf k}=dk_xdk_ydk_z$ standing for volume element.

\subsection{Representation by a superposition of monochromatic cylindrical waves -- by the Fourier--Bessel integral}

Now we pass in \eqref{uk} to cylindrical coordinates,
  $x=\rho\cos\varphi\,, y=\rho\sin\varphi\,,$
$x\cos\phi+y\sin\phi=\rho \cos(\phi-\varphi)\,.$
Integration with respect to   $\phi$ and application of the well-known expression for the Bessel function    $J_0$,  $J_0(m)=\frac{1}{2\pi}\int_0^{2\pi}e^{im\cos{\mu}}d\mu$, see \cite{Abram}, provides
\begin{gather*}\label{}
u = -i\int_{0}^{\frac{\pi}{2}} {\sin{\cal{X}}\,d{\cal{X}}} \int_{0}^{\infty} \hat{f}(k \cos{\cal{X}}) J_0(k\rho\sin{\cal{X}})
e^{ik[(ct+ib) -(z+ib)\cos{\cal{X}}]}k\,dk \,.
\end{gather*}
By substitution
$  k_z=k\cos{\cal{X}}$
we present \eqref{RNDW} as follows
\begin{gather}\label{cyl}
u=-i\int_{0}^{\infty}e^{ik(ct+ib)} dk \int_{0}^{k} \, \hat{f}(k_z)
J_0(\sqrt{k^2-k_z^2}\,\rho)e^{-ik_z(z+ib)}dk_z\,.
\end{gather}
Introducing the variable $\omega=ck$ we come up with
\begin{gather}
u=-\frac{i}{c}\int_{0}^{\infty}e^{i\omega (t+ib/c)} d\omega \int_{0}^{\omega/ c} \, \hat{f} \left(k_z\right)
J_0(\sqrt{(\omega/c)^2-k_z^2}\,\,\rho)e^{-ik_z(z+ib)}dk_z\,
\end{gather}
which allows the following relation between the waveform $f(\theta)$ in \eqref{RNDW} and the weight $A(k_z,\omega)$ in \eqref{FouBess}:
 \begin{equation}A(k_z,\omega)=-\frac{i}{c}e^{-(\omega /c-k_z)b}   \, \hat{f}(k_z).
 \end{equation}

\section{Conclusions}

Thus, we have established relationships between several
representations of localized unidirectional waves. These are   quasi-spherical waves \eqref{RNDW},  Fourier-Bessel integrals \eqref{FouBess} and  superpositions of monochromatic \eqref{Fourier3D} and non-stationary \eqref{uF1} plane waves.

\section*{Acknowledgments}
The authors are grateful to I.\,Besieris, P.\,Saari, and J.\,Lekner
  for a useful exchange of views, and to
  N.N.\,Rosanov, and I.A.\,So  for valuable comments.

\section*{Conflict of interest}
The authors declare that they have no conflict of interest.
{

\end{document}